\documentclass[11pt]{article}
\usepackage{amssymb,latexsym}
\usepackage{epsfig}
\usepackage{eufrak}
\usepackage{amsmath}
\usepackage{mathrsfs}
\usepackage{color}

\newtheorem{theorem}{Theorem}
\newtheorem{lemma}{Lemma}

\newcommand{\be}{\begin{equation}}
\newcommand{\ee}{\end{equation}}
\newcommand{\bee}{\begin{eqnarray*}}
\newcommand{\eee}{\end{eqnarray*}}
\newcommand{\bel}{\begin{eqnarray}}
\newcommand{\eel}{\end{eqnarray}}
\newcommand{\bec}{\begin{cases}}
\newcommand{\eec}{\end{cases}}
\newcommand{\bem}{\begin{bmatrix}}
\newcommand{\eem}{\end{bmatrix}}

\newcommand{\la}{\label}
\newcommand{\li}{\left}
\newcommand{\ri}{\right}

\newcommand{\vol}{\mathrm{vol}}
\newcommand{\ep}{\epsilon}

\newcommand{\lm}{\lambda}

\newcommand{\de}{\delta}
\newcommand{\De}{\Delta}
\newcommand{\vDe}{\varDelta}

\newcommand{\se}{\theta}

\newcommand{\al}{\alpha}

\newcommand{\ro}{\rho}
\newcommand{\ka}{\kappa}
\newcommand{\om}{\omega}
\newcommand{\Om}{\Omega}

\newcommand{\f}{\frac}

\newcommand{\cd}{\cdots}

\newcommand{\qu}{\quad}
\newcommand{\qqu}{\qquad}
\newcommand{\fa}{\forall}

\newcommand{\mscr}{\mathscr}
\newcommand{\mcal}{\mathcal}
\newcommand{\mbf}{\mathbf}
\newcommand{\bb}{\mathbb}

\newcommand{\mrm}{\mathrm}
\newcommand{\bs}{\boldsymbol}

\newcommand{\ap}{\approx}

\newcommand{\tx}{\text}

\newcommand{\iy}{\infty}

\newcommand{\pa}{\partial}

\newcommand{\bed}{\begin{description}}
\newcommand{\eed}{\end{description}}
\newcommand{\bei}{\begin{itemize}}
\newcommand{\eei}{\end{itemize}}
\newcommand{\ben}{\begin{enumerate}}
\newcommand{\een}{\end{enumerate}}
\newcommand{\bib}{\bibitem}
\newcommand{\beL}{\begin{lemma}}
\newcommand{\eeL}{\end{lemma}}
\newcommand{\beT}{\begin{theorem}}
\newcommand{\eeT}{\end{theorem}}
\newcommand{\sect}{\section}

\newcommand{\bpf}{\begin{pf}}
\newcommand{\epf}{\end{pf}}
\newcommand{\bsk}{\bigskip}
\newcommand{\bi}{\binom}

\newcommand{\pfbox}{\hfill\mbox{$\Box$}}

\newenvironment{pf}{\paragraph*{Proof{\rm.}}}{\pfbox\bigskip}

\begin{document}

\title{{\bf Robust Estimation of Mean Values}
\thanks{The author had been previously working with Louisiana State University at Baton Rouge, LA 70803, USA,
and is now with Department of Electrical Engineering, Southern
University and A\&M College, Baton Rouge, LA 70813, USA; Email:
chenxinjia@gmail.com}}

\author{Xinjia Chen}

\date{November 2008}

\maketitle

\begin{abstract}

In this paper, we develop a computational approach for estimating
the mean value of a quantity in the presence of uncertainty. We
demonstrate that, under some mild assumptions, the upper and lower
bounds of the mean value are efficiently computable via a sample
reuse technique, of which the computational complexity is shown to
posses a Poisson distribution.

\end{abstract}

\sect{Introduction}

In many situations, it is desirable to estimate the mean value of a
scalar quantity $\bs{Q}$ which is a function of independent random
vectors $\bs{V}$ and $\bs{\vDe}$ such that the distribution of
$\bs{V}$ is known and that the distribution of $\bs{\vDe}$ is
unknown \cite{Huber}. Namely, it is interested to estimate the
expectation of $\bs{Q} = q(\bs{V}, \bs{\vDe})$, where $q(.,.)$ is a
multivariate function. From modeling considerations, it is
reasonable to assume that $\bs{\vDe}$ is bounded in norm $||.||$,
and radially symmetrical and nondecreasing in its probability
density function, $f_{ \bs{\vDe} } (.)$ with the following notions:

(i) The norm, $||\bs{\vDe}||$, of $\bs{\vDe}$ is no greater than a
certain value $r$, i.e., $||\bs{\vDe}|| \leq r$;

(ii) For any realization $\vDe$ of $\bs{\vDe}$, $f_{ \bs{\vDe} }
(\vDe)$ depends only on, $||\vDe||$, the norm of $\vDe$;

(iii) For any $\vDe_1$ and $\vDe_2$ such that $||\vDe_1|| <
||\vDe_2||$, $f_{ \bs{\vDe} } (\vDe_1) \geq f_{ \bs{\vDe} }
(\vDe_2)$.

Such assumptions have been proposed by Barmish and Lagoa
\cite{Barmish} in the context of robustness analysis of control
systems, where $\bs{\vDe}$ is referred to as ``uncertainty'' because
of the lack of knowledge of its distribution.

In this paper, we shall focus on the estimation of the expectation
$\bb{E} [ \bs{Q} ] = \bb{E} [ q(\bs{V}, \bs{\vDe}) ]$ based on
assumptions (i), (ii) and (iii).  Such a problem is referred to as
{\it robust estimation} due to the fact that the exact distribution
of $\bs{\vDe}$ is not available. In the special case that the
maximum norm $r$ of $\bs{\vDe}$ equals $0$, the robust estimation
problem reduces to a conventional estimation problem.  Instead of
seeking the exact value of $\bb{E} [ \bs{Q} ]$ which is obviously
impossible, we aim at obtaining upper and lower bounds for $\bb{E} [
\bs{Q} ]$. It is intuitive that the gap between the upper and lower
bounds should be increasing with respect to $r$. Since the relation
between $Q$ and $\bs{V}, \; \bs{\vDe}$ can be fairly complicated,
the Monte Carlo estimation method is the unique and powerful
approach.

The remainder of the paper is organized as follows.  In Section 2,
we derive upper and lower bounds for $\bb{E} [ \bs{Q} ]$ based on
assumptions (i), (ii) and (iii).  In Section 3, we propose a Monte
Carlo method for the evaluation of the bounds of $\bb{E} [ \bs{Q}
]$. In particular, we introduce a sample reuse method to
substantially reduce the computational complexity.  In Section 4, we
investigate the computational complexity of the Monte Carlo method
implemented with the principle of sample reuse.  Section 5 is the
conclusion.

\sect{Bounds of Expectation}

In this section, we shall derive upper and lower bounds of $\bb{E} [
\bs{Q} ] = \bb{E} [q(\bs{V}, \bs{\vDe})]$ based on the assumptions
described in Section 1.  For this purpose, we have the following
fundamental result, which is a slight generalization of the uniform
principle proposed by Barmish and Lagoa \cite{Barmish}.

\beT  Let $\bs{\vDe}_\ro^{\mrm{u}}$ be a random vector with a
uniform distribution over $\{ \vDe: || \vDe || \leq \ro \}$.  Define
\[
\bb{M} (\ro) = \bb{E} \li [ q (\bs{V}, \bs{\vDe}_\ro^{\mrm{u}}) \ri
], \qu \underline{\bb{M}} (r) = \inf_{0 < \ro < r} \bb{M} (\ro), \qu
\overline{\bb{M}} (r) = \sup_{0 < \ro < r} \bb{M} (\ro).
\]
Then, $\underline{\bb{M}} (r) < \bb{E} [ \bs{Q} ] <
\overline{\bb{M}} (r)$.  \eeT

\bsk

See Appendix A for a proof. Theorem 1 reveals that the computation
of the bounds of $\bb{E} [ \bs{Q} ]$ can be reduced to the
evaluation of function $\bb{M} (\ro)$, which can be accomplished via
Monte Carlo simulation.  A conventional  method is as follows:

Partition interval $(0, r]$ by grid points $r = \ro_1 > \ro_2 > \cd
> \ro_m > 0$.  Let $\mscr{B}_\ell = \{ \vDe: ||\vDe|| \leq \ro_\ell
\}$.  For $\ell = 1, \cd, m$, estimate $\bb{M} (\ro_\ell)$ as the
empirical mean \[ \f{ \sum_{i = 1}^N q (\bs{V}_i, X_{\ell, i} ) } {
N }
\]
where $\bs{V}_i, \; X_{\ell, i}, \; i = 1, \cd, N$ are mutually
independent random variables such that $\bs{V}_1, \cd, \bs{V}_N$ are
i.i.d. random samples of $\bs{V}$ and $ X_{\ell, 1}, \cd, X_{\ell,
N}$ are i.i.d. random samples uniformly distributed over
$\mscr{B}_\ell$.  Clearly, the total number of simulations is $Nm$
for estimating $\bb{M} (\ro_\ell), \; \ell = 1, \cd, m$.  A major
problem with this approach is that the computational complexity can
be extremely high, since the number of grid points $m$ is typically
a very large number. To overcome such a problem, we shall develop a
sample reuse technique in the next section.

\sect{Sample Reuse}

In this section, we shall explore the idea of sample reuse to reduce
the computational complexity.  The sample reuse method has been
proposed by Chen et al. \cite{C0, C1} for the robustness analysis of
control systems.  The idea of sample reuse is to start simulation
from the largest set $\mscr{B}_1$ and if it also belongs to smaller
subsets the experimental result is saved for later use in the
smaller sets. As can be seen from last section, a conventional
approach would require a total of $N m$ simulations. However, due to
sample reuse, the actual number of experiments for set
$\mscr{B}_\ell$ is a random number $\mbf{n}_\ell$, which is usually
much less than $N$. Hence, this strategy saves a significant amount
of computational effort.

In order to provide a precise description of the principle of sample
reuse, we assume that all random variables are defined in the same
probability space $(\Om, \mscr{F}, \Pr )$. We shall introduce a
function $\mscr{G}$, referred to as {\it sample reuse function}, as
follows.

Let $X_1, \cd, X_m$ be i.i.d. samples uniformly distributed over
$\mscr{A}$.  Let $Y_1, \cd, Y_n$ be i.i.d. samples uniformly
distributed over $\mscr{B}$. Let $m \leq n $ and $\mscr{A} \supset
\mscr{B}$. Define reusable sample size $\mbf{k}$ such that $\mbf{k}
(\om)$ is the number of elements of $\{ X_i(\om) \in \mscr{B}:  i =
1, \cd, m \}$ for any $\om \in \Om$. Define random variables $Z_1,
\cd, Z_n$ such that, for any $\om \in \Om$,  \[ Z_\ell (\om) = \bec
X_{i_\ell} (\om) & \tx{for} \; 1 \leq \ell \leq \mbf{k} (\om),\\
Y_\ell (\om) & \tx{for} \; \mbf{k} (\om) < \ell \leq n \eec
\]
where $i_\ell, \; 1 \leq \ell \leq \mbf{k} (\om)$ are the indexes of
the elements of $\{ X_i(\om) \in \mscr{B}: i = 1, \cd, n \}$ such
that $i_\ell$ is increasing with respect to $\ell$. This process of
generating $Z_1, \cd, Z_n$ from $X_1, \cd, X_m$ and $Y_1, \cd, Y_n$
is denoted by
\[
(Z_1, \cd, Z_n; \mbf{k}) = \mscr{G} (X_1, \cd, X_m; Y_1, \cd, Y_n).
\]
With regard to the distribution of $Z_1, \cd, Z_n$, we have \beT
Suppose $X_1, \cd, X_m$ are independent with $Y_1, \cd, Y_n$. Then,
$Z_1, \cd, Z_n$ are i.i.d. samples uniformly distributed over
$\mscr{B}$. \eeT

See Appendix B for a proof. Now we can use $\mscr{G}$ to precisely
describe the sample reuse algorithm for estimating $\bb{M}
(\ro_\ell), \; \ell = 1, \cd, m$. Let $X_{\ell,i}, \; i = 1, \cd, N$
be the random samples uniformly distributed over $\mscr{B}_\ell$ for
$\ell = 1, \cd, m$.  Let $Y_{1,i} = X_{1,i}$ for $i = 1, \cd, N$ and
$(Y_{\ell,1}, \cd, Y_{\ell, N}; \mbf{k}_\ell) = \mscr{G} (
Y_{\ell-1, 1}, \cd, Y_{\ell-1, N}; X_{\ell, 1}, \cd, X_{\ell, N} )$
for $\ell = 2, \cd, m$. As a result of Theorem 1, we have that, for
any $\ell \in \{1, \cd, m\}$, random variables $Y_{\ell, i}, \; i =
1, \cd, N$ have the same associated cumulative distribution with
that of random variables $X_{\ell, i}, \; i = 1, \cd, N$.  This
implies that $\f{1}{N} \sum_{i = 1}^N q (\bs{V}_i, Y_{\ell, i} )$
has the same distribution as that of $\f{1}{N} \sum_{i = 1}^N q
(\bs{V}_i, X_{\ell, i} )$ for $\ell = 1, \cd, m$.  Therefore, we can
use $\f{1}{N} \sum_{i = 1}^N q (\bs{V}_i, Y_{\ell, i} )$ as an
estimator of $\bb{M} (\ro_\ell)$ for $\ell = 1, \cd, m$.  By virtue
of such sample reuse method, the total number of simulations is
reduced from $Nm$ to $N + \sum_{\ell = 2}^m \mbf{n}_\ell$, where
$\mbf{n}_\ell = N - \mbf{k}_\ell$ for $\ell = 2, \cd, m$. As will be
demonstrated in the next section, this can be a huge reduction of
complexity for a large $m$.

\sect{Poisson Complexity}

Since the total number of simulations for using the sample reuse
method to estimate $\bb{M} (\ro_\ell), \; \ell = 1, \cd, m$ is $N +
\sum_{\ell = 2}^m \mbf{n}_\ell$, it is important to investigate the
distribution of $\sum_{\ell = 2}^m \mbf{n}_\ell$. In this regard, we
have the following general result.

\beT \la{basic} For arbitrary sequence of nested sets $\mscr{B}_1
\supset \mscr{B}_2 \supset \cd \supset \mscr{B}_m$ with
 $\mrm{vol} (\mscr{B}_1) = V_{\mrm{max}}$ and $\mrm{vol} (\mscr{B}_m)
= V_{\mrm{min}}$, the cumulative distribution function of
$\sum_{\ell = 2}^m \mbf{n}_\ell$ is bounded from below by the
cumulative distribution function of a Poisson random variable
$\bs{P}$ with mean {\small $\lm = N \ln \li ( \f{ V_{\mrm{max}} } {
V_{\mrm{min}} } \ri )$}. That is, {\small $\Pr \li \{ \sum_{\ell =
2}^m \mbf{n}_\ell = 0 \ri \} = \Pr \{ \bs{P} = 0 \}$} and {\small
$\Pr \li \{ \sum_{\ell = 2}^m \mbf{n}_\ell \leq k \ri \} > \Pr \{
\bs{P} \leq k \}$} for any positive integer $k$. Moreover, as the
maximum difference of volumes of all consecutive sets tends to be
zero, $\sum_{\ell = 2}^m \mbf{n}_\ell$ converges to $\bs{P}$ in
distribution. \eeT

\bsk

See Appendix C for a proof. It should be noted that the volume of a
set $\mscr{B}$, denoted by $\mrm{vol} ( \mscr{B})$, is referred to
the Lebesgue measure of $\mscr{B}$ in this paper.

As an immediate consequence of Theorem 3, we have
\[ \Pr \li \{ \sum_{\ell = 2}^m \mbf{n}_\ell  > 0 \ri \} = \Pr \{
\bs{P} > 0 \}, \qqu \Pr \li \{ \sum_{\ell = 2}^m \mbf{n}_\ell > k
\ri \} < \Pr \{ \bs{P} > k \}, \qu k = 1, 2, \cd
\]
which implies that
\[
\bb{E} \li [ \sum_{\ell = 2}^m \mbf{n}_\ell \ri ] = \sum_{k = 0}^\iy
\Pr \li \{ \sum_{\ell = 2}^m \mbf{n}_\ell
> k \ri \} < \sum_{k = 0}^\iy \Pr \{ \bs{P} > k \} = \lm = N \ln \li (
\f{ V_{\mrm{max}} } { V_{\mrm{min}} } \ri ).
\]
By virtue of Theorem 3, we can derive some simple bounds for the
distribution of $\sum_{\ell = 2}^m \mbf{n}_\ell$ as follows. \beT
$\Pr \{ \sum_{\ell = 2}^m \mbf{n}_\ell \geq k \} \leq e^{- \lm} \li
( \f{ \lm e } { k } \ri )^k$ for any number {\small $k > \lm = N \ln
\li ( \f{ V_{\mrm{max}} } { V_{\mrm{min}} } \ri )$}. In particular,
$\Pr \{ \sum_{\ell = 2}^m \mbf{n}_\ell \geq e \lm \} \leq e^{- \lm}$
and $\Pr \{ X \geq (1 + \ep) \lm \} < \exp \li ( - \f{ \ep^2 \lm } {
4 } \ri )$ for $0 < \ep < 1$.
 \eeT

\bsk

See Appendix D for a proof.

Now we apply Theorem 3 to investigate the density of original
samples of $\bs{\vDe}$.   Suppose that the volume of $\{ \vDe: ||
\vDe || \leq \ro \}$ is proportional to $\ro^d$ where $d$ is the
dimension of the set.  Let $\bs{N}_\ro$ denote the number of
original samples included in $\{ \vDe: || \vDe || \leq \ro \}$ when
applying the sample reuse method to interval $[ \f{\ro}{\ka}, \ro
]$. Define the density of samples at radius $\ro$ as $\mcal{D} (\ro)
= \lim_{\de \to 0} \f{ \bb{E} [\bs{N}_{\ro + \de} -
\bs{N}_\ro]}{\de}$.  Then, we have the following result.

\beT $\mcal{D}(\ro)$ is equal to $\f{N d}{\ro} \li ( \f{\ka \ro}{a}
\ri)^d$ for $\ro \in (0, \f{a}{\ka} ]$ and is less than $\f{N
d}{\ro}$ for $\ro \in (\f{a}{\ka}, a]$. \eeT

\bsk

\begin{figure}
\centering
\includegraphics[height=10cm]{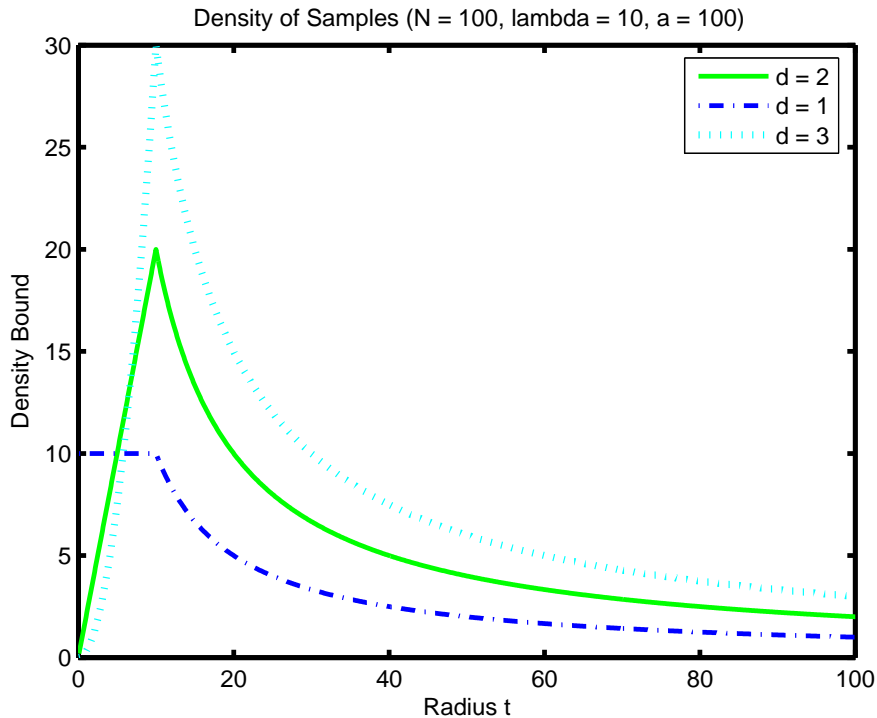}
\caption{Illustrative Example ($N = 100,  \; a = 100, \; \lm = 10$)}
\label{density}
\end{figure}

See Appendix E for a proof. From this theorem, we can obtain an
upper bound for the expected number of original samples with norm
bounded in $[0, a]$.  As can be seen from Theorem 5, the density
function is  unimodal and achieves the largest value at $\ro =
\f{a}{\ka}$.  The density function is displayed by Figure
\ref{density}.

\sect{Conclusion}

We have proposed an efficient computational approach for estimating
the mean value of a random function, for which the distribution of
relevant random variables are not completely available.  A Monte
Carlo method with sample reuse as a key mechanism is established.
The associated computational complexity is demonstrated to follow a
Poisson distribution.

\bsk

\appendix

\sect{Proof of Theorem 1}

We follow the similar method of Barmish and Lagoa \cite{Barmish}.
Let $\mcal{V}$ denote the volume of $\mcal{B} = \{ \vDe: || \vDe ||
\leq r \}$. We partition the set $\mcal{B}$ as $K$ layers of equal
volume $\f{\mcal{V}}{K}$ such that the $k$-th layer is $\mcal{L}_k =
\{ \vDe: r_{k - 1} < || \vDe || \leq r_k \}$ with $0 = r_0 < r_1 <
r_2 < \cd < r_K = r$. Then, the density function can be expressed as
 \[
f_{\bs{\vDe}} (\De ) \ap \sum_{k =1}^K \bb{I}_k(\De) \lm_k
 \]
 where $\lm_k, \; k = 1, \cd, K$ satisfying
 \be
 \la{cont}
\f{\mcal{V}}{K} \sum_{k = 1}^K \lm_k = 1, \qqu \lm_1 \geq \lm_2 \geq
\cd \geq \lm_K \geq 0
 \ee
 and $\bb{I}_k(.)$ is the indicator function such that $\bb{I}_k(\De) = 1$ if $\De$ falls into the $k$-th
 layer $\mcal{L}_k$ and $\bb{I}_k(\De) = 0$ otherwise.  Let $f_{\bs{V} } (.)$ denote the density function of $\bs{V}$.
Since $\bs{V}$ and $\bs{\vDe}$ are independent, we have
 \bee \bb{E} [q(\bs{V}, \bs{\vDe}) ]
& \ap & \int_{\{(v, \De): ||\De|| \leq r \}} q(v, \De) \;  f_{\bs{V} } (v) dv \; f_{\bs{\vDe}} (\De) d \De\\
& = & \int_{\{(v, \De): ||\De|| \leq r \}} q(v, \De) \;  f_{\bs{V} }
(v) dv \; \li [ \sum_{k =1}^K \bb{I}_k(\De) \lm_k \ri ] d \De =
\sum_{k =1}^K  \al_k  \lm_k,
 \eee
where $\al_k =  \int_{\{(v, \De): ||\De|| \leq r \}} q(v, \De) \;
\bb{I}_k(\De) f_{\bs{V} } (v) dv \;
 d \De$.  Therefore, the upper and lower bounds of $\bb{E} [q(\bs{V}, \bs{\vDe})
 ]$ correspond to the maximum and minimum of the linear program: $\sum_{k =1}^K  \al_k
 \lm_k$ subject to constraint (\ref{cont}).  From convex analysis,
 the maximum and minimum of this linear program are achieving at extreme points of the form:
 \[
\lm_k = \bec \f{K}{j \mcal{V}} & \tx{for} \; 1 \leq k \leq j,\\
0 & \tx{for} \; j < k \leq K. \eec
 \]
As the number of layers $K$ tends to infinity, the summation
$\sum_{k =1}^K \bb{I}_k(\De) \lm_k$, which is associated with
extreme point $(\lm_1, \cd, \lm_K)$,  tends to a uniform
distribution. This justifies the theorem.

\sect{Proof of Theorem 2}

Let $S_\ell \subseteq {\mscr{B}}$ for $\ell = 1, \cd, n$. Define
$\bb{D} = \{1, \cd, n\}$ and $\mcal{I}_s = \{ (i_1, \cd, i_s): i_1 <
\cd < i_s; \; i_\ell \in \bb{D}, \; \ell = 1, \cd, s \}$. Then,
{\small \bee \Pr \{ Z_\ell \in S_\ell, \; \ell = 1, \cd, n \} & = &
\sum_{s = 0}^n \sum_{ (i_1, \cd, i_s) \in \mcal{I}_s } \Pr \{
X_{i_\ell} \in S_\ell, \; \ell = 1, \cd, s; \; X_j \notin
{\mscr{B}}, \; j
\in \bb{D} \setminus \{i_1, \cd, i_s\} \}\\
&  & \times \Pr \{ Y_\ell \in S_\ell, \; \ell = s + 1, \cd, n \}.
\eee} For simplicity of notations, we let $V_{S_\ell} = \mrm{vol} (
S_\ell ), \; V_{\mscr{A}} =  \mrm{vol} ( \mscr{A} )$ and
$V_{\mscr{B}} =  \mrm{vol} ( \mscr{B} )$.  Note that $\Pr \{ Y_\ell
\in S_\ell, \; \ell = s + 1, \cd, n \} = \prod_{\ell = s + 1}^n \li
( \f{V_{S_\ell}}{V_{\mscr{B}}} \ri )$ and {\small \bee \Pr \{
X_{i_\ell} \in S_\ell, \; \ell = 1, \cd, s; \; X_j \notin
{\mscr{B}}, \; j \in \bb{D} \setminus \{i_1, \cd, i_s\} \} & = & \li
( \f{V_{\mscr{A}} - V_{\mscr{B}}}{V_{\mscr{A}}} \ri )^{m - s}
\prod_{\ell = 1}^s \li (
\f{V_{S_\ell}}{V_{\mscr{A}}} \ri )\\
& = & \li ( \f{ V_{\mscr{B}} } { V_{\mscr{A}} } \ri )^s \li ( 1 -
\f{ V_{\mscr{B}}}{V_{\mscr{A}}} \ri )^{m - s} \prod_{\ell = 1}^s \li
( \f{V_{S_\ell}}{V_{\mscr{B}}} \ri ). \eee} Since there are
$\bi{n}{s}$ elements in $\mcal{I}_s$, we have {\small \bee \Pr \{
Z_\ell \in S_\ell, \; \ell = 1, \cd, n \} & = & \sum_{s = 0}^n
\bi{n}{s} \li ( \f{ V_{\mscr{B}} } { V_{\mscr{A}} } \ri )^s \li ( 1
- \f{ V_{\mscr{B}}}{V_{\mscr{A}}} \ri )^{m - s} \prod_{\ell = 1}^s
\li ( \f{V_{S_\ell}}{V_{\mscr{B}}} \ri ) \prod_{\ell = s + 1}^n \li
(
\f{V_{S_\ell}}{V_{\mscr{B}}} \ri )\\
& = & \prod_{\ell = 1}^n \li ( \f{V_{S_\ell}}{V_{\mscr{B}}} \ri )
\sum_{s = 0}^n \bi{n}{s} \li ( \f{ V_{\mscr{B}} } { V_{\mscr{A}} }
\ri )^s \li ( 1 - \f{
V_{\mscr{B}}}{V_{\mscr{A}}} \ri )^{m - s}\\
& = & \prod_{\ell = 1}^n \li ( \f{V_{S_\ell}}{V_{\mscr{B}}} \ri ).
\eee} This concludes the proof of the theorem.

\sect{Proof of Theorem 3}

We need some preliminary results.

\beL \la{indep} Let $N_1 \leq N_2 \leq \cd \leq N_m$. For $\ell = 1,
\cd, m$, let $v_\ell = \mrm{vol} ( \mscr{B}_\ell )$ and $X_{\ell,i},
\; i = 1, \cd, N_\ell$ be i.i.d. random samples uniformly
distributed over $\mscr{B}_\ell$. Let $Y_{1,i} = X_{1,i}$ for $i =
1, \cd, N_1$ and $(Y_{\ell,1}, \cd, Y_{\ell, N_\ell}; \mbf{k}_\ell)
= \mscr{G} ( Y_{\ell-1, 1}, \cd, Y_{\ell-1, N_{\ell- 1}}; X_{\ell,
1}, \cd, X_{\ell, N_\ell} )$ for $\ell = 2, \cd, m$. Define
$\mbf{n}_\ell = N_\ell - \mbf{k}_\ell$ for $\ell = 2, \cd, m$.
 Then, {\small $\Pr
\{ \mbf{n}_\ell = n_\ell, \; \ell = 2, \cd, m \} = \prod_{\ell =
2}^m B \li ( N_{\ell} - n_{\ell}, N_{\ell -1}, \f{v_\ell} {v_{\ell -
1}} \ri )$} for $N_{\ell} - N_{\ell - 1} \leq n_\ell \leq N_\ell$
and $2 \leq \ell \leq m$, where $B(k, n, p) = \bi{n}{k} p^k (1 -
p)^{n - k}$. \eeL

\bpf  We use induction method.  First, it is easy to show that the
lemma is true for $m =2$. Next, we assume that the lemma is true for
$m-1$ and show that the lemma is also true for $m$. Let $\Pr \{(k_1,
\cd, k_m), \; (N_1, \cd, N_m), \; (v_1, \cd, v_m) \}$ denote the
probability that, among the $N_1$ samples generated from the biggest
set $\mscr{B}_1$,
 there are
$k_\ell$ samples falling into $\mscr{B}_\ell$ for $\ell = 1, 2, \cd,
m$.  Let $\mbf{P}^m \li \{ (n_2, \cd, n_m), \; (N_1, \cd, N_m), \;
(v_1, \cd, v_m) \ri \}$ denote the probability of event
$\{\mbf{n}_\ell = n_\ell, \; \ell = 2, \cd, m \}$ associated with
the application of the sample reuse method to sets $\mscr{B}_\ell,
\; \ell = 1, \cd, m$ with required sample sizes $N_1 \leq N_2 \leq
\cd \leq N_m$.  Let $\mbf{P}^{m - 1} \{ (n_3, \cd, n_m), \; (N_2 -
k_2, \cd, N_m - k_m), \; (v_2, \cd, v_m)  \}$ denote the probability
of event $\{\mbf{n}_\ell = n_\ell, \; \ell = 3, \cd, m \}$
associated with the application of the sample reuse method to sets
$\mscr{B}_\ell, \; \ell = 2, \cd, m$ with required sample sizes $N_2
- k_2 \leq  \cd \leq N_m - k_m$.  Note that \bee & & \mbf{P}^m \li
\{ (n_2, \cd, n_m), \; (N_1, \cd, N_m), \; (v_1, \cd, v_m) \ri \}\\
& = & \sum _{k_2 \geq
k_3 \geq \cd \geq k_m \geq 0} \Pr \{(k_1, \cd, k_m), \;  (N_1, \cd, N_m), \; (v_1, \cd, v_m)  \}\\
&    & \times \; \mbf{P}^{m-1} \li \{ (n_3, \cd, n_m), \; \li (N_2 -
k_2, \; N_3 - k_3, \; \cd, \; N_m - k_m \ri ), \; (v_2, \cd, v_m)
\ri \}
 \eee
where $n_2 + k_2 = N_2$ and $k_1 = N_1$.  By the mechanism of sample
reuse, {\small \[ \Pr \{(k_1, \cd, k_m), \; (N_1,
\cd, N_m), \; (v_1, \cd, v_m)  \}\\
  =  \li [ \prod_{\ell = 2}^{m}
\bi{k_{\ell - 1}}{k_{\ell}} \li ( \f{v_{\ell - 1} - v_{\ell}}{v_1}
\ri )^{k_{\ell - 1} - k_{\ell} }  \ri ] \li ( \f{ v_m } { v_1 } \ri
)^{k_m}. \]} Since $N_\ell$ and $- k_\ell$ are non-decreasing with
respect to $\ell$, we have that $N_\ell - k_\ell$ is non-decreasing
with respect to $\ell$. Hence, by the assumption of induction,
{\small \bee &  & \mbf{P}^{m-1} \li \{ (n_3, \cd, n_m), \; \li (N_2
- k_2, \; N_3 - k_3, \; \cd, \; N_m - k_m \ri ), \; (v_2, \cd, v_m)
\ri \}\\
& = & \prod_{\ell = 3}^m B \li ( N_{\ell} - n_{\ell} -
k_{\ell}, N_{\ell -1} - k_{\ell - 1}, \f{v_\ell} {v_{\ell
- 1}} \ri )\\
& = &  \prod_{\ell = 3}^m \bi{N_{\ell-1} - k_{\ell-1}}{ N_\ell -
n_\ell - k_\ell} \li ( \f{ v_\ell } { v_{\ell-1} } \ri )^{N_\ell -
n_\ell - k_\ell} \li ( 1 - \f{ v_{\ell} } { v_{\ell-1} } \ri
)^{N_{\ell-1} - N_\ell + n_{\ell} - k_{\ell-1} + k_\ell} \eee} and
consequently, {\small \bee & & \mbf{P}^m \li \{ (n_2, \cd, n_m),
\; (N_1, \cd, N_m), \; (v_1, \cd, v_m) \ri \}\\
& = &  \sum _{k_2 \geq k_3 \geq \cd \geq k_m \geq 0} \li [
\prod_{\ell = 2}^{m} \bi{k_{\ell - 1}}{k_{\ell}} \li ( \f{v_{\ell -
1} - v_{\ell}}{v_1} \ri )^{k_{\ell - 1} - k_{\ell} }  \ri ] \li (
\f{ v_m } { v_1 } \ri )^{k_m}\\
&    & \times  \prod_{\ell = 3}^m \bi{N_{\ell-1} - k_{\ell-1}}{
N_\ell - n_\ell - k_\ell} \li ( \f{ v_\ell } { v_{\ell-1} } \ri
)^{N_\ell - n_\ell - k_\ell} \li ( 1 - \f{ v_{\ell} } { v_{\ell-1} }
\ri )^{N_{\ell-1} - N_\ell + n_{\ell} - k_{\ell-1} + k_\ell}\\
& = & \sum _{k_2 \geq k_3 \geq \cd \geq k_m \geq 0} \li [
\prod_{\ell = 2}^{m} \bi{k_{\ell - 1}}{k_{\ell}} \bi{N_{\ell-1} -
k_{\ell-1}}{ N_\ell - n_\ell -
k_\ell} \ri ] \\
&   & \times \li ( \f{ v_m } { v_1 } \ri )^{k_m} \prod_{\ell =
2}^{m} \li ( \f{v_{\ell - 1} - v_{\ell}}{v_1} \ri )^{k_{\ell - 1} -
k_{\ell} } \times \prod_{\ell = 3}^m \li ( \f{ v_\ell } { v_{\ell-1}
} \ri )^{N_\ell - n_\ell - k_\ell} \li ( 1 - \f{ v_{\ell} } {
v_{\ell-1} } \ri )^{N_{\ell-1} - N_\ell + n_{\ell} - k_{\ell-1} +
k_\ell}.  \eee} Making use of the relationships $k_1 = N_1$ and $k_2
= N_2 - n_ 2$, we have {\small \bee &   &  \li ( \f{ v_m } { v_1 }
\ri )^{k_m} \prod_{\ell = 2}^{m} \li ( \f{v_{\ell - 1} -
v_{\ell}}{v_1} \ri )^{k_{\ell - 1} - k_{\ell} }  \times \prod_{\ell
= 3}^m \li ( \f{ v_\ell } { v_{\ell-1} } \ri )^{ - k_\ell} \li ( 1 -
\f{ v_{\ell} } { v_{\ell-1} } \ri )^{ -
k_{\ell-1} + k_\ell}\\
& = & \li ( v_{1} - v_{2} \ri )^{k_{1} - k_{2} } \li ( \f{ v_2^{k_2}
} { v_1^{N_1} } \ri ) = \li ( \f{ v_2 } { v_{1} } \ri )^{N_2 - n_2}
\li ( 1 - \f{ v_{2} } { v_{1} } \ri )^{N_{1} - N_2 + n_{2}} \eee}
and thus {\small \bee &   & \li ( \f{ v_m } { v_1 } \ri )^{k_m}
\prod_{\ell = 2}^{m} \li ( \f{v_{\ell - 1} - v_{\ell}}{v_1} \ri
)^{k_{\ell - 1} - k_{\ell} } \times \prod_{\ell = 3}^m \li ( \f{
v_\ell } { v_{\ell-1} } \ri )^{N_\ell - n_\ell - k_\ell} \li ( 1 -
\f{ v_{\ell} } { v_{\ell-1} } \ri )^{N_{\ell-1} - N_\ell +
n_{\ell} - k_{\ell-1} + k_\ell}\\
& = & \li ( \f{ v_m } { v_1 } \ri )^{k_m} \prod_{\ell = 2}^{m}  \li
( \f{v_{\ell - 1} - v_{\ell}}{v_1} \ri )^{k_{\ell - 1} - k_{\ell} }
 \times \prod_{\ell = 3}^m \li ( \f{ v_\ell } { v_{\ell-1} } \ri )^{
- k_\ell} \li ( 1 - \f{ v_{\ell} } { v_{\ell-1} } \ri )^{ -
k_{\ell-1} + k_\ell}\\
&    & \times \prod_{\ell = 3}^m \li ( \f{ v_\ell } { v_{\ell-1} }
\ri )^{N_\ell - n_\ell} \li ( 1 - \f{ v_{\ell} } { v_{\ell-1} } \ri
)^{N_{\ell-1} - N_\ell + n_{\ell}}\\
& = & \prod_{\ell = 2}^m \li ( \f{ v_\ell } { v_{\ell-1} } \ri
)^{N_\ell - n_\ell} \li ( 1 - \f{ v_{\ell} } { v_{\ell-1} } \ri
)^{N_{\ell-1} - N_\ell + n_{\ell}}. \eee} On the other hand, if the
lemma really holds, we have {\small \bee & & \mbf{P}^{m} \li \{
(n_2, \cd, n_m), \; (N_1, \cd, N_m), \; (v_1,
\cd, v_m) \ri \}\\
& = & \prod_{\ell = 2}^m B \li ( N_{\ell} - n_{\ell}, N_{\ell -1},
\f{v_\ell} {v_{\ell - 1}} \ri ) =  \li [ \prod_{\ell = 2}^m
\bi{N_{\ell-1} }{ N_\ell - n_\ell} \ri ]  \li [ \prod_{\ell = 2}^m
\li ( \f{ v_\ell } { v_{\ell-1} } \ri )^{N_\ell - n_\ell} \li ( 1 -
\f{ v_{\ell} } { v_{\ell-1} } \ri )^{N_{\ell-1} - N_\ell + n_{\ell}}
\ri ].  \eee} Therefore, to show the lemma, it remains to show
{\small
\[ \sum_{k_2 \geq
k_3 \geq \cd \geq k_m \geq 0} \li [ \prod_{\ell = 2}^m \bi{k_{\ell -
1}}{k_\ell} \bi{N_{\ell-1} - k_{\ell-1}}{ N_\ell - n_\ell - k_\ell}
\ri ]
 =  \prod_{\ell = 2}^m \bi{N_{\ell-1}}{ N_\ell - n_\ell}. \]}
Using the relationships $k_1 = N_1$ and $k_2 = N_2 - n_ 2$,  this
identity can be reduced to  the following identity \[ \sum_{k_2 \geq
k_3 \geq \cd \geq k_m \geq 0} \li [ \prod_{\ell = 3}^m \bi{k_{\ell -
1}}{k_\ell} \bi{N_{\ell-1} - k_{\ell-1}}{ N_\ell - n_\ell - k_\ell}
\ri ]
 =  \prod_{\ell = 3}^m \bi{N_{\ell-1}}{ N_\ell - n_\ell}, \]
 which can be shown by observing that
{\small \bee &  & \sum_{k_2 \geq k_3 \geq \cd \geq k_{m - i} \geq 0}
\li [ \prod_{\ell = 3}^{m - i} \bi{k_{\ell - 1}}{k_\ell}
\bi{N_{\ell-1} -
k_{\ell-1}}{ N_\ell - n_\ell - k_\ell} \ri ]\\
 & = & \sum_{k_2 \geq
k_3 \geq \cd \geq k_{{m - i}-1} \geq 0} \li [ \prod_{\ell = 3}^{{m -
i}-1} \bi{k_{\ell - 1}}{k_\ell} \bi{N_{\ell-1} - k_{\ell-1}}{ N_\ell
- n_\ell - k_\ell} \ri ] \sum_{k_{m - i} = 0}^{k_{{m - i}-1}}
\bi{k_{{m - i} - 1}}{k_{m - i}}
\bi{N_{{m - i}-1} - k_{{m - i}-1}}{ N_{m - i} - n_{m - i} - k_{m - i}}\\
& = & \sum_{k_2 \geq k_3 \geq \cd \geq k_{{m - i}-1} \geq 0} \li [
\prod_{\ell = 3}^{{m - i}-1} \bi{k_{\ell - 1}}{k_\ell}
\bi{N_{\ell-1} - k_{\ell-1}}{ N_\ell - n_\ell - k_\ell} \ri ]
\bi{N_{{m - i}-1}}{ N_{m - i} - n_{m - i}} \eee} for $0 \leq i \leq
m - 4$ and {\small \[ \sum_{k_2 \geq k_3 \geq 0} \li [ \prod_{\ell =
3}^{3} \bi{k_{\ell - 1}}{k_\ell} \bi{N_{\ell-1} - k_{\ell-1}}{
N_\ell - n_\ell - k_\ell} \ri ] = \sum_{k_2 \geq k_3 \geq 0} \li [
\bi{k_{2}}{k_3} \bi{N_{2} - k_{2}}{ N_3 - n_3 - k_3} \ri ] =
\bi{N_{2}}{ N_3 - n_3}.  \]} This completes the proof of the lemma.
\epf

\beL \la{bet} Let $\se > 1$ and $N \geq 1$. Define $L(\se, k) =
\sum_{i = 0}^k \bi{N}{i} \li ( 1 - \f{1}{\se} \ri )^i \li (
\f{1}{\se} \ri )^{N- i}$ and $L_P(\se, k) = \sum_{i=0}^k \f{(N \ln
\se)^i } { i! }  \exp(- N \ln \se )$ for $k = 0, 1, \cd, N$.   Then,
$L(\se, 0) = L_P(\se, 0)$ and $L(\se, k) > L_P(\se, k)$ for $k = 1,
\cd, N$.  \eeL

\bpf First, it is evident that $L( \se, 0) = L_P(\se, 0) = \se^{-N}$
and $L(\se, N) = 1 > L_P(\se, N)$. Hence, it remains to show the
lemma for $k = 1, \cd, N -1$. It is easy to show that $\lim_{\se \to
\iy} L(\se, k) = \lim_{\se \to \iy} L_P(\se, k) = 0$ and thus
$\lim_{\se \to \iy} [L(\se, k) - L_P(\se, k)] = 0$ for $k = 1, \cd,
N -1$.  It can also be readily checked that $\lim_{\se \to 1} L(\se,
k) = \lim_{\se \to 1} L_P(\se, k) = 1$ and consequently $\lim_{\se
\to 1} [L(\se, k) - L_P(\se, k)] = 0$ for $k = 1, \cd, N -1$. Noting
that {\small $\f{ \pa L (\se, k) } { \pa \se  }  = - \f{ N! } {
k!(N-k - 1)! } \li ( 1 - \f{1}{\se} \ri )^k \li ( \f{1}{\se} \ri
)^{N-k + 1}$} and {\small $\f{ \pa L_P (\se, k) } { \pa \se  }  = -
\f{ ( N \ln \se )^k } { k! } \f{N}{\se^{N+1}}$},  we have {\small
$\f{ \pa [ L (\se, k) - L_P(\se, k)]} { \pa \se  } = \f{N} {k!
\se^{N+1}} \li [ (N \ln \se )^k - \f{ (N -1) ! (\se - 1)^k } {(N - k
- 1)! } \ri ] > 0$} if and only if $\varphi(\se) > 0$, where
$\varphi(\se) = \ln \se - \al (\se -1)$ with {\small $\al = \li [
\f{(N-1)!  } { (N - k - 1)! } \ri ]^{\f{1}{k}} \f{1}{N} < 1$}. Since
$\varphi(1) = 0$ and $\f{d \varphi(\se)}{d \se} = \f{1}{\se} - \al$
is positive for $\se \in \li ( 1, \f{1}{\al} \ri )$,  we have
$\varphi(\se) > 0$ for $\se \in \li ( 1, \f{1}{\al} \ri ]$. Since
$\varphi \li ( \f{1}{\al} \ri ) > 0, \; \lim_{\se \to \iy}
\varphi(\se) < 0$ and $\f{d \varphi(\se)}{d \se} < 0$ for $\se >
\f{1}{\al}$,  there exists a unique number $\se^*$ greater than
$\f{1}{\al}$ such that $\varphi(\se^*) = 0$. Hence, $\varphi(\se)$
is positive for $\se \in (1, \se^*)$ and negative for $\se > \se^*$.
This implies that $L (\se, k) - L_P(\se, k)$ is monotonically
increasing with respect to $\se \in (1, \se^*)$ and monotonically
decreasing with respect to $\se \in ( \se^*, \iy)$. Recalling that
$\lim_{\se \to 1} [ L (\se, k) - L_P(\se, k) ] = \lim_{\se \to \iy}
[ L (\se, k) - L_P(\se, k) ] = 0$, we have $L (\se, k)
> L_P(\se, k)$ for any $\se > 1$.  This completes the proof of the
lemma.

 \epf

\beL \la{ADD} Let $U_i,\; V_i, \; i = 1, \cd, n$ be mutually
independent non-negative discrete random variables. Suppose that
$\Pr \{U_i = 0 \} = \Pr \{ V_i = 0 \}$ and $\Pr \{U_i \leq k\} > \Pr
\{ V_i \leq k \}$ for any positive integer $k$ and $i = 1, \cd, n$.
Then, $\Pr \{ \sum_{i = 1}^n U_i = 0 \} = \Pr \{ \sum_{i = 1}^n V_i
= 0 \}$ and $\Pr \{ \sum_{i = 1}^n U_i \leq k \}
> \Pr \{ \sum_{i = 1}^n V_i \leq k \}$ for any positive integer $k$.
 \eeL

 \bpf  We use induction method.  The lemma is obviously true for $n = 1$.
 Assuming that the lemma is true for $n = m - 1 \geq 1$, we have
 {\small $\Pr \{ \sum_{i =
 1}^{m} U_i = 0 \}  =  \Pr \{ \sum_{i =
 1}^{m - 1} U_i = 0, \; U_m  = 0 \} =    \Pr \{ \sum_{i =
 1}^{m - 1} V_i = 0 \} \Pr \{V_m = 0 \} =  \Pr \{ \sum_{i =
 1}^{m} V_i = 0 \}$} and
 {\small $\Pr \{ \sum_{i =
 1}^{m} U_i \leq k \} = \sum_{l = 0}^k \Pr \{ \sum_{i =
 1}^{m - 1} U_i = l, \; U_m \leq k - l \}
  >  \sum_{l = 0}^k \Pr \{ \sum_{i =
 1}^{m - 1} V_i = l  \} \Pr \{V_m \leq k - l \} =  \Pr  \{ \sum_{i =
 1}^{m} V_i \leq k  \}$ } for any positive integer $k$, which implies that the lemma is also
true for $n = m$.  By the principle of induction, the lemma is
established.

 \epf

We are now in a position to prove the theorem.  We shall first show
that the distribution of $\sum_{\ell = 2}^m \mbf{n}_\ell$ is bounded
from below by the distribution of a Poisson variable with mean
{\small $N \ln \f{ V_{\mrm{max}} } { V_{\mrm{min}} }$}.  Define $U_i
= \mbf{n}_{i + 1}$ for $i = 1, \cd, m - 1$. Then, by Lemma
\ref{indep}, $U_i$ are independent binomial random variables such
that $\Pr \{ U_i \leq k \} = L(\se_i, k)$ for $k = 0, 1, \cd, N$ and
$i = 1, \cd, m - 1$, where $\se_i = \f{v_i}{v_{i + 1}}$. Define
Poisson variables $V_i, \; i = 1, \cd, m - 1$ such that $U_i, \;
V_i, \; i = 1, \cd, m - 1$ are mutually independent and that $\Pr \{
V_i \leq k \} = L_P(\se_i, k)$ for non-negative integer $k$ and $i =
1, \cd, m - 1$. By Lemmas \ref{bet} and \ref{ADD}, we have {\small
$\Pr \{ \sum_{\ell = 2}^m \mbf{n}_\ell = 0 \} = \Pr \{ \sum_{i =
1}^{m - 1} U_i = 0 \} = \Pr \{ \sum_{i = 1}^{m - 1} V_i = 0 \}$} and
{\small $\Pr \{ \sum_{\ell = 2}^m \mbf{n}_\ell \leq k \} = \Pr \{
\sum_{i = 1}^{m - 1} U_i \leq k \} > \Pr \{ \sum_{i = 1}^{m - 1} V_i
\leq k \}$} for $k = 1, 2, \cd$. Noting that $V_1, \cd, V_{m - 1}$
are independent Poisson variables with corresponding means $N \ln
\se_1, \cd, N \ln \se_{m - 1}$,  we have that $\sum_{i = 1}^{m - 1}
V_i$ is also a Poisson variable with mean $N \sum_{i = 1}^{m - 1}
\ln \se_i = N \sum_{i = 1}^{m - 1} \ln \f{v_i}{v_{i + 1}} = N \ln
\f{v_1}{v_m} = N \ln \f{ V_{\mrm{max}} } { V_{\mrm{min}} }$.

Next, we shall show that the distribution of $\sum_{\ell = 2}^m
\mbf{n}_\ell$ tends to be the distribution of a Poisson variable
with mean {\small $N \ln \f{ V_{\mrm{max}} } { V_{\mrm{min}} }$} as
$\nu = \max \{ v_\ell - v_{\ell + 1}: 1 \leq \ell \leq m - 1 \}$,
the maximum difference between the volumes of two consecutive nested
sets, tends to be zero while the volumes of $\mscr{B}_1$ and
$\mscr{B}_m$ respectively assume fixed values $v_1 = V_{\mrm{max}}$
and $v_m = V_{\mrm{min}}$.

Since all sample sizes are equal to $N$, by Lemma \ref{indep}, for
$\ell = 2, \cd, m$, the original sample sizes $\mbf{n}_\ell, \; \ell
= 2, \cd, m$ are mutually independent binomial random variables such
that $\Pr \{ \mbf{n}_\ell = k  \} = B(k, N, p_\ell)$ for $0 \leq k
\leq N$ and $2 \leq \ell \leq m$, where {\small $p_\ell = 1 -
\f{v_\ell } {v_{\ell - 1}}$} with $v_\ell = \vol (\mscr{B}_\ell)$.
Therefore, the moment generating function of $\sum_{\ell =2}^m
\mbf{n}_\ell$ can be expressed as $G (s) = \li [ \prod_{\ell = 2}^m
(p_\ell s + 1 - p_\ell) \ri ]^N$, where $s \in (0, 1]$ is a real
number.  Since $p_\ell s + 1 - p_\ell$ is positive for any $s \in
(0,1]$ and $\ell = 2, \cd, m$, it is meaningful to define $g(s) =
\sum_{\ell = 2}^m \ln (p_\ell s + 1 - p_\ell)$ for $s \in (0,1]$.
Hence, $G(s) = \exp( N g(s) )$.  For simplicity of notations, define
{\small $h(s) = (s - 1) \ln \li ( \f{ V_{\tx{max}} } { V_{\tx{min}}
} \ri ), \; I_1(s) = \int_{0}^s \sum_{\ell = 2}^m \f{ v_{\ell - 1} -
v_\ell } { z ( v_{\ell - 1} - v_\ell ) + v_\ell} d z - \int_{0}^1
\sum_{\ell = 2}^m \f{ v_{\ell - 1} - v_\ell } { z ( v_{\ell - 1} -
v_\ell ) + v_\ell} d z$} and {\small $I_2 (s) = \int_{0}^s
\sum_{\ell = 2}^m  \f{ v_{\ell - 1} - v_\ell } { v_\ell} d z -
\int_{0}^1 \sum_{\ell = 2}^m  \f{ v_{\ell - 1} - v_\ell } { v_\ell}
d z$}.  The lemma can be established by the following three steps.

First, it can be seen that $g(s) = I_1(s)$ for any $s \in (0,1]$,
since $I_1(1) = g(1) = 0$ and {\small \bee \f{d I_1(s) } { d s } =
\sum_{\ell = 2}^m  \f{ v_{\ell - 1} - v_\ell } { s ( v_{\ell - 1} -
v_\ell ) + v_\ell} = \sum_{\ell = 2}^m \f{p_\ell} {p_\ell s + 1 -
p_\ell } = \f{d g(s) } { d s }\eee} for any $s \in (0,1]$.

Second, we need to show that $| I_1 (s ) - I_2 (s) | \to 0$ for any
$s \in (0,1]$ as $\nu \to 0$.  Noting that {\small \bee \li |
\int_0^s \f{ z (v_{\ell - 1} - v_\ell)^2 } {  v_\ell^2 + z ( v_{\ell
- 1} - v_\ell ) v_\ell } d z \ri |  & = & \f{ (v_{\ell - 1} -
v_\ell)^2 } { v_\ell }  \li | \int_0^s \f{ z } {  v_\ell + z (
v_{\ell - 1} - v_\ell )
 } d z \ri |\\
 & \leq & \f{ (v_{\ell - 1} - v_\ell)^2 } { v_\ell }  \int_0^s
\f{ z } {  v_\ell
 } d z =   \f{ s^2 (v_{\ell - 1} - v_\ell)^2 } { 2 v_\ell^2 } \leq
 \f{ s^2 \nu (v_{\ell - 1} - v_\ell) } { 2 V_{\mrm{min}}^2 }
\eee} for any $s \in (0,1]$, we have {\small \bee | I_1 (s ) - I_2
(s) | & \leq & \sum_{\ell = 2}^m  \li | \int_0^s \f{ z (v_{\ell - 1}
- v_\ell)^2 } { v_\ell^2 + z ( v_{\ell - 1} - v_\ell ) v_\ell } d z
\ri | + \sum_{\ell = 2}^m \li | \int_0^1 \f{ z (v_{\ell - 1} -
v_\ell)^2
} { v_\ell^2 + z ( v_{\ell - 1} - v_\ell ) v_\ell } d z \ri |\\
& \leq & \sum_{\ell = 2}^m  \f{ s^2 \nu (v_{\ell - 1} - v_\ell) } {
2 V_{\mrm{min}}^2 } +  \sum_{\ell = 2}^m \f{ \nu (v_{\ell - 1} -
v_\ell) } { 2 V_{\mrm{min}}^2 }\\
& = & \f{ (s^2 + 1) \nu  } { 2 V_{\mrm{min}}^2 } \sum_{\ell = 2}^m
(v_{\ell - 1} - v_\ell)  = \f{ (s^2 + 1) ( V_{\mrm{max}} -
V_{\mrm{min}} ) \nu  } {2 V_{\mrm{min}}^2 }.\eee} Therefore, $| I_1
(s ) - I_2 (s) | \to 0$ for any $s \in (0,1]$ and arbitrary $v_\ell,
\; \ell = 1, \cd, m$, as $\nu \to 0$.

Third, we need to show $g(s) \to h(s)$ as $\nu \to 0$.  Since
{\small \bee h(s) - I_2 (s) & = & \int_{0}^s \li [ \ln \li ( \f{
V_{\tx{max}} } { V_{\tx{min}} } \ri ) - \sum_{\ell = 2}^m \f{
v_{\ell - 1} - v_\ell } {
 v_\ell}  \ri ] d z  - \int_{0}^1 \li [  \ln \li ( \f{
V_{\tx{max}} } { V_{\tx{min}} } \ri ) - \sum_{\ell = 2}^m \f{
v_{\ell - 1} - v_\ell } { v_\ell} \ri ] d z, \eee} we have {\small
$| I_2 (s) - h(s) |  \leq  \int_{0}^s \li | \ln \li ( \f{
V_{\tx{max}} } { V_{\tx{min}} } \ri ) - \sum_{\ell = 2}^m \f{
v_{\ell - 1} - v_\ell } { v_\ell}  \ri | d z  + \int_{0}^1 \li | \ln
\li ( \f{ V_{\tx{max}} } { V_{\tx{min}} } \ri ) - \sum_{\ell = 2}^m
\f{ v_{\ell - 1} - v_\ell } { v_\ell} \ri | d z$}.  By the
definition of Riemann integration, {\small $\sum_{\ell = 2}^m \f{
v_{\ell - 1} - v_\ell } { v_\ell} \to \int_{V_{\mrm{min}}}^{
V_{\mrm{max}} } \f{d v}{v}  = \ln \li ( \f{ V_{\tx{max}} } {
V_{\tx{min}} } \ri )$} as $\nu \to 0$ for arbitrary $v_\ell, \; \ell
= 1, \cd, m$. It follows that, for any $s \in (0,1]$ and arbitrary
$v_\ell, \; \ell = 1, \cd, m$, $| I_2 (s ) - h (s) | \to 0$ as $\nu
\to 0$.  In view of $| g(s) - h(s) | = | I_2(s) - h(s) + I_1(s) -
I_2(s) | \leq | I_2(s) - h(s) | + | I_1(s) - I_2(s) |$, we have
$g(s) \to h(s)$ as $\nu \to 0$ for any $s \in (0,1]$ and arbitrary
$v_\ell, \; \ell = 1, \cd, m$.  Therefore, we can conclude that
{\small $G(s)  \to \exp \li( (s-1) N \ln \li ( \f{ V_{\tx{max}} } {
V_{\tx{min}} } \ri ) \ri )$} as $\nu \to 0$ for any $s \in (0,1]$
and arbitrary $v_\ell, \; \ell = 1, \cd, m$.  This proves that
 $\sum_{\ell = 2}^m \mbf{n}_\ell$ converges in distribution to a
Poisson variable of mean {\small $N \ln \li ( \f{ V_{\tx{max}} } {
V_{\tx{min}} } \ri )$}. The proof of the theorem is thus completed.

\sect{Proof of Theorem 4}

We need some preliminary results.

\beL Let $X$ be a Poisson variable of mean $\lm > 0$. For any number
$k > \lm$, $\Pr \{ X \geq k \} \leq e^{- \lm} \li ( \f{ \lm e } { k
} \ri )^k$.
 \eeL

\bpf Since $\Pr \{ X \geq k \} = \Pr \li \{ e^{ t (X - k)  } \geq 1
\ri \} \leq \bb{E} \li [ e^{ t (X - k)  }  \ri ]$ for any $t > 0$,
we have
 $\Pr \{ X \geq k \} \leq \inf_{t > 0} \bb{E} \li [ e^{ t (X - k) }
\ri ]$.  Note that \bee \bb{E} \li [ e^{ t (X - k) }  \ri ] =
\sum_{i = 0}^\iy e^{ t (i - k)} \f{ \lm^i } { i! } e^{- \lm} =
e^{\lm e^t } e^{- \lm} e^{ - t k } \sum_{i = 0}^\iy \f{ (\lm e^t)^i
} { i! } e^{- \lm e^t} =  e^{- \lm} e^{ \lm e^t - t k }, \eee which
is minimized if and only if $\lm e^t = k$.  Since $k > \lm$, we have
$t = \ln \li ( \f{k}{\lm} \ri ) > 0$ such that $\lm e^t = k$.  For
this value of $t$, we have $e^{- \lm} e^{ \lm e^t - t k } = e^{-
\lm} \li ( \f{ \lm e } { k } \ri )^k$.  Hence, we have shown $\Pr \{
X \geq k \} \leq e^{- \lm} \li ( \f{ \lm e } { k } \ri )^k$.

\epf

\bsk

Now we are in a position to prove the theorem.  By Theorem 3, we
have $\Pr \{ \sum_{\ell = 2}^m \mbf{n}_\ell \geq k \} \leq \Pr \{ X
\geq k \} \leq e^{- \lm} \li ( \f{ \lm e } { k } \ri )^k$. Setting
$k = e \lm$, we have $\Pr \{ \sum_{\ell = 2}^m \mbf{n}_\ell \geq e
\lm \} \leq e^{- \lm}$.  Moreover, using the inequality $(1 + \ep)
\ln (1 + \ep) > \ep + \f{ \ep^2 } {4}, \; \fa \ep \in (0,1]$, we
have $\Pr \{ \sum_{\ell = 2}^m \mbf{n}_\ell \geq (1 + \ep) \lm \} <
\li [  \f{ e^\ep } { (1 + \ep)^{1 + \ep} }  \ri ]^\lm < \exp \li ( -
\f{ \ep^2 \lm } { 4 } \ri )$ for $0 < \ep < 1$.  This completes the
proof of the theorem.

\sect{Proof of Theorem 5} By Theorem 3, we have {\small $\bb{E}
[\bs{N}_{\ro + \de}] \leq N \li [ 1 + d \ln \li ( \f{ \ka (\ro +
\de)  } { a }\ri ) \ri ]$}. Now fix the gridding over $(\ro, \ro +
\de]$. By Theorem 3,  as the griding over $[\f{a}{\ka}, \ro]$
becomes increasingly dense, we have {\small $\bb{E} [\bs{N}_{\ro}]
\to N \li [ 1 + d \ln \li ( \f{ \ka \ro  } { a }\ri ) \ri ]$}. This
implies that, for any $\ep
>0$, we have {\small $\bb{E} [\bs{N}_{\ro}] > N \li [ 1 + \ln \li (
\f{ \ka \ro } { a }\ri ) \ri ] - \ep$} for a sufficiently dense
gridding over $[\f{a}{\ka}, \ro]$. Hence, {\small  \bee \bb{E}
[\bs{N}_{\ro + \de} - \bs{N}_\ro] = \bb{E} [\bs{N}_{\ro + \de}] -
\bb{E} [\bs{N}_\ro]
 <  N d \ln \li ( \f{ \ka (\ro + \de)  } { a }\ri ) - N d \ln \li ( \f{ \ka \ro } { a }\ri ) + \ep
  =  N d \ln \li ( \f{\ro + \de}{\ro} \ri ) + \ep. \eee} Since the
argument holds for any small $\ep >0$, we have {\small $\bb{E}
[\bs{N}_{\ro + \de} - \bs{N}_\ro] \leq N d \ln \li ( \f{\ro +
\de}{\ro} \ri )$}.  Therefore, the density {\small $\mcal{D}(\ro) =
\lim_{\de \to 0} \f{ \bb{E} [\bs{N}_{\ro + \de} - \bs{N}_\ro] }
{\de} \leq \lim_{\de \to 0} \f{ N d \ln \li ( \f{\ro + \de}{\ro} \ri
) } { \de } = \f{N d}{\ro}$}. On the other hand, as the gridding
gets dense, we have {\small $\bb{E} [\bs{N}_{\ro + \de} -
\bs{N}_\ro] \to N d \ln \li ( \f{\ro + \de}{\ro} \ri )$} and thus
$\mcal{D}(\ro) \to \f{N d}{\ro}$.  For $\ro \in (0, \f{a}{\ka} ]$,
it follows from Theorem 3 that $\bs{N}_\ro$ is a binomial random
variable corresponding to $N$ i.i.d. trials with a success
probability $\li( \f{\ka \ro}{a} \ri)^d$.  Hence, $\bb{E}
[\bs{N}_\ro] = N \li ( \f{\ka \ro}{a} \ri)^d$ and accordingly
{\small $\mcal{D}(\ro) = \f{N d}{\ro} \li ( \f{\ka \ro}{a} \ri)^d$}.
This completes the proof of the theorem.

\end{document}